\documentclass[11pt]{article}
\usepackage{graphicx}

\usepackage[centertags]{amsmath}
\usepackage{amsfonts}
\usepackage{amssymb}
\usepackage{newlfont}
\usepackage{enumerate}
\usepackage[ansinew]{inputenc}

\newtheorem{theorem}{Theorem}[section]

\newtheorem{lemma}[theorem]{Lemma}

\newtheorem{problem}[theorem]{Problem}

\newtheorem{conjecture}[theorem]{Conjecture}

\def   \proof          {\noindent{\it Proof}. }

\def   \qed            {\hfill $\Box$\par\bigskip}
\def\Box{\hbox{\vrule width6pt height6pt depth0pt}}

\newcommand{\Nset}{{\mathbb N}}
\newcommand{\Rset}{{\mathbb R}}

\newcommand{\comment}[1]{}

\begin{document}

\title{Complements of nearly perfect graphs}

\author{
\normalsize Andr\'as Gy\'arf\'as\thanks{Research supported in part by OTKA Grant No. K104343}\\
\small Alfr\'ed R\'enyi Institute of Mathematics\\[-0.8ex]
\small Hungarian Academy of Sciences\\[-0.8ex]
\small Budapest, P.O. Box 127, H-1364, Hungary\\[-0.8ex]
\small \texttt{gyarfas.andras@renyi.mta.hu}\\[-0.8ex]\\
\and Zhentao Li \thanks{This research partially supported by the FQRNT. } \\
\small \'Ecole Normale Sup\'erieure de Lyon, LIP , \'Equipe MC2
\\[-0.8ex]
\small 46, all\'ee d`Italie, 69342 Lyon Cedex 07, France\\[-0.8ex]
\small {\tt zhentao.li@ens-lyon.fr}\\
\and Raphael Machado \\
\small {Inmetro, Rio de Janeiro, Brazil}\\
\small {\tt rcmachado@inmetro.gov.br}\\
\and Andr\'as Seb\H o\\
\small CNRS, Grenoble-INP, UJF, Laboratoire G-SCOP, \\[-0.8ex]
%\small CNRS\\[-0.8ex]
\small 46 avenue F\'elix Viallet, 38031 Grenoble Cedex, France \\[-0.8ex]
\small \texttt{andras.sebo@g-scop.inpg.fr}\\
\and St\'ephan Thomass\'e \\
\small \'Ecole Normale Sup\'erieure de Lyon, LIP , \'Equipe MC2
\\[-0.8ex]
\small 46, all\'ee d`Italie, 69342 Lyon Cedex 07, France\\[-0.8ex]
\small {\tt stephan.thomasse@ens-lyon.fr}\\
\and Nicolas Trotignon\thanks{partially supported by the French Agence Nationale de la
      Recherche under reference \textsc{anr-10-jcjc-Heredia}} \\
\small CNRS, \'Ecole Normale Sup\'erieure de Lyon, LIP , \'Equipe MC2
\\[-0.8ex]
\small 46, all\'ee d`Italie, 69342 Lyon Cedex 07, France\\[-0.8ex]
\small \texttt{nicolas.trotignon@ens-lyon.fr}
}

%\date{}
\maketitle
\newpage

\begin{abstract}
  A class of graphs closed under taking induced subgraphs is
  $\chi$-bounded if there exists a function $f$ such that for all
  graphs $G$ in the class, $\chi(G) \leq f(\omega(G))$.  We consider
  the following question initially studied in [A.~Gy{\'a}rf{\'a}s,
  \newblock Problems from the world surrounding perfect graphs, {\em
    Zastowania Matematyki Applicationes Mathematicae}, 19:413--441,
  1987].  For a $\chi$-bounded class $\cal C$, is the class $\overline{C}$
  $\chi$-bounded (where $\overline{\cal C}$ is the class of graphs formed
  by the complements of graphs from $\cal C$)?

  We show that if $\cal C$ is $\chi$-bounded by the constant function
  $f(x)=3$, then $\overline{\cal C}$ is $\chi$-bounded by
  $g(x)=\lfloor\frac{8}{5}x\rfloor$ and this is best possible.  We
  show that for every constant $c>0$, if $\cal C$ is $\chi$-bounded by a
  function $f$ such that $f(x)=x$ for $x \geq c$, then $\overline{\cal C}$
  is $\chi$-bounded. For every $j$, we construct a class of
  graphs $\chi$-bounded by $f(x)=x+x/\log^j(x)$ whose complement is not $\chi$-bounded.
\end{abstract}

\section{Introduction}
\label{s:introduction}

In the present paper, we consider simple and finite graphs.  We denote
by $\chi(G)$ (resp. $\theta(G)$) the minumum number of stable sets
(resp.\ cliques) needed to cover the vertices of $G$.  We denote by
$\omega(G)$ (resp.\ $\alpha(G)$) the maximum size of a clique (resp.\
stable set) of $G$.  A graph $G$ is \emph{$\chi$-bounded} (resp.\
\emph{$\theta$-bounded}) by a function $f$ if $\chi(H)\leq
f(\omega(H))$ (resp.\ $\theta(H)\leq f(\alpha(H))$) for every induced
subgraph $H$ of $G$.  Observe that a graph $G$ is $\chi$-bounded by
$f$ if an only if its complement $\overline{G}$ is $\theta$-bounded by
$f$.  A class of graphs is \emph{$\chi$-bounded} (resp.\
\emph{$\theta$-bounded}) if for some function $f$, every graph of the
class is $\chi$-bounded (resp.\ $\theta$-bounded) by $f$.  The class
of \emph{perfect graphs} is the class of graph $\chi$-bounded by the
identity function.  For every $\chi$-bounded (resp.\ $\theta$-bounded)
class, there exists a smallest $\chi$-bounding (resp.\
$\theta$-bounding) function that we refer as the \emph{optimal
  $\chi$-bounding (resp.\ $\theta$-bounding) function} for the class.

We address a general question asked by Gy{\'a}rf{\'a}s~\cite{Gyarfas}:
for which functions $f$ is the class of graphs $\chi$-bounded by $f$ also
$\theta$-bounded (by a possibly different function $g$)?  Such
functions are called \emph{complementary-bounded} functions and $g$ is
a \emph{complementary bounding} function for $f$.  If $f$ is a
complementary bounding function, we denote by $f^*$ the optimal
$\theta$-bounding function of the class of graphs $\chi$-bounded by
$f$.

\begin{theorem}[K\H onig \cite{konig:31}]
  \label{th:k}
  If $G$ is bipartite, then $\theta(G) = \alpha(G)$.
\end{theorem}

The classical theorem above can be rephrased as ``the identity is the
optimal complementary bounding function for the constant function
$f=2$'', or by $2^* = \text{id}$.  In Section~\ref{sect:3col}, we push
further this line of research by computing the optimal complementary
bounding function of the constant function $f=3$.  To do so, we prove
that 3-colourable graphs are $\theta$-bounded by $f^*(x)=
\left\lfloor\frac{8}{5}x\right\rfloor$ or equivalently $3^*=\lfloor {8\over 5}
\text{id} \rfloor$ (Theorem~\ref{thm:3col}).  Our proof uses a
well-known result of Gallai~\cite{gallai:colorCritical} on color
critical graphs (Theorem~\ref{th:Gcritcial} below).  This result has
been sharpened by Stehl\'\i k~\cite{stehlik:03}.  Using ideas of
Cornu\'ejols, Hartvigsen and Pulleyblank \cite{CoP:83} the sharpened
theorem can be closely related to matching theory and namely to
Gallai~\cite{gallai:factorCritical}.  These relations, a short proof
using Theorem~\ref{th:gallai}, and a systematic account of
reformulations and generalizations of these results are being laid
out by Seb\H{o} and Stehl\'\i k \cite{sebostehlik:hypergraphs}.

The following remarks on the constant function $f_m=m$ are
from~\cite{Gyarfas}.  On one hand, $f_m^*(x)\le \lfloor{m+1\over
  2}\rfloor x$ because the vertex set of any $m$-chromatic graph can
be covered by the vertices of at most $\lfloor{m+1\over 2}\rfloor$
bipartite graphs. On the other hand, $f_m^*(x)\ge {mx\over 2}$ for
$x>x_0(m)$ from a nice probabilistic construction of Erd\H os
\cite{ER}: for arbitrary $m$, there exists $t$ and an $m$-partite
triangle-free graph $G$ with $t$ vertices in each partite class, such
that $\alpha(G)=t$.  In fact, our construction on Fig.~\ref{fig:g58}
showing $f_3^*(x)\ge \left\lfloor\frac{8}{5}x\right\rfloor$ is such a
graph with $m=3,t=5$.  Thus we have the asymptotic of $f_m^*$ for even
$m$'s, but not for odd $m$'s, apart from $m=3$.  It seems hard to make
an intelligent guess even on the asymptotic $f_5^*(x)$.

\medskip

K\H onig's theorem was generalized by Lov\'asz in the Perfect Graph
Theorem~\cite{Lovasz}, stating that the identity is its own optimal
complementary bounding function.  A function $f$ is \emph{eventually
  identity} if there exist a constant $c$ such that for all $x \ge c$,
$f(x) = x$.  In Section~\ref{sec:evc}, we prove eventually
identity functions are complementary-bounded (Theorem~\ref{old}).
This theorem was stated without proof in~\cite{Gyarfas}.  Our proof is
an induction that reduces the problem to Lov\'asz's theorem.

\medskip

In Section~\ref{sec:nontb}, we deal with functions that are not
complementary-bounded.  Gy{\'a}rf{\'a}s~\cite{Gyarfas} proved that for
every real number $\varepsilon>0$, the function $f(x) = x +
\varepsilon x$ (and thus any function greater than $f$) is not
complementary bounded.  We improve this result by proving that the
function $f(x)=x+x/\log^j(x)$ is
not-complementary-bounded for any $j$ (Theorem~\ref{th:unbound2}). The methods we use
to prove $\chi$-boundedness in this section rely again on a theorem of
Gallai~\cite{gallai:factorCritical} on {\em factor-critical graphs}, graphs where removing any vertex yields a graph with a perfect matching.

\medskip

It was conjectured in~\cite{Gyarfas} that $f(x)=x+c$ is complementary-bounded for any constant $c$.  This conjecture remains open
even for $c=1$.  Our tools from Sections~\ref{sect:3col}
and~\ref{sec:evc} are not strong enough to prove that $f$ is
complementary-bounded, and our tools from Section~\ref{sec:nontb} are
not strong enough to prove that it is not.

\section{Complementary bounding function of $3$-chromatic
  graphs}
\label{sect:3col}

In this section, we find the smallest $\theta$-bounding function for
the class of 3-colourable graphs.  We work on a more general class
$\cal C$: graphs $G$ such that for every induced subgraph $H$ of $G$,
$\alpha(H) \geq |V(H)|/3$.  Graphs satisfying this property are more
general than 3-colourable graphs as the 5-wheel has this property but is not 3-colourable.

\begin{theorem}
  \label{thm:3col}
  Every graph $G$ in $\cal C$ satisfies $\theta(G) \le \lfloor
  \frac{8}{5}\alpha(G) \rfloor$.  This is best possible, in the sense
  that for every integer $x\geq 0$, there exists a graph $G$ in $\cal
  C$ with $\alpha(G) = x$ and $\theta(G) = \lfloor
  \frac{8}{5}\alpha(G) \rfloor$.
\end{theorem}

\begin{figure}
  \begin{center}
  \includegraphics{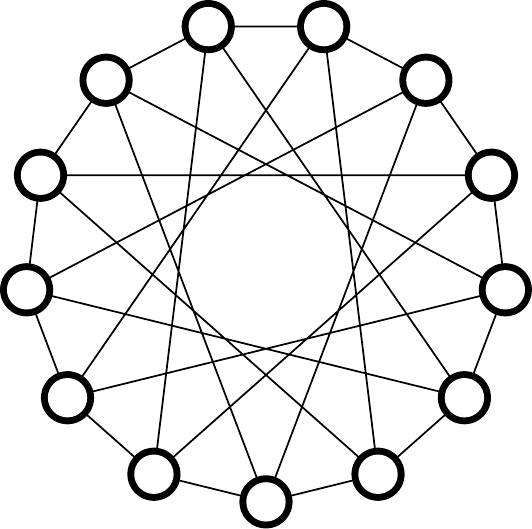}
  \rule{2em}{0ex}
  \includegraphics{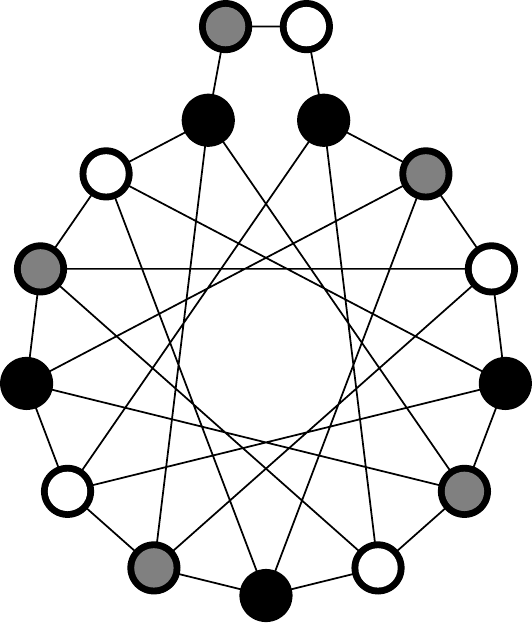}
  \end{center}
  \caption{Graphs $R_{3,5}$ and  $G_{5,8}$}
  \label{fig:g58}
\end{figure}

Our result improves on the previous upper bound of $\frac{5}{3}x$
from~\cite{Gyarfas}.  The rest of the section is devoted to proving
Theorem~\ref{thm:3col}.  It is best possible because of the graph
$G_{5, 8}$ represented in Fig.~\ref{fig:g58} satisfies
$$|V(G_{5,8})|=15,\  \omega(G_{5,8})=2,\ \chi(G_{5,8})=3,\
\alpha(G_{5,8})=5,\ \theta(G_{5,8})=8.$$
Other graphs with the same parameters can also be found as induced subgraphs in some of the seven graphs with parameters $|V(G)|=17, \omega(G)=2, \alpha(G)=5$ (one is
given in~\cite{kery:ramsey}, all seven in~\cite{kalbfleisch:edge}).
However, $G_{5,8}$ is much simpler for our purposes.  Checking
$|V(G_{5,8})|=15$, $\chi(G_{5,8})=3$ and $\theta(G_{5,8})\leq 8$ is
immediate from the figure (cycles of length~5 are easy to find, a
3-colouring is shown and a clique cover of size 8 is obtained by
taking every second edge on the obvious hamiltonian cycle and an
isolated vertex).  Note that $\chi(G_{5,8}) = 3$ implies that $G\in
{\cal C}$.  To compute $\alpha$ and $\omega$, it is convenient to
consider the graph $R_{3, 5}$, also represented in Fig.~\ref{fig:g58},
that is well known in Ramsey Theory as the unique graph $G$ on at
least 13 vertices such that $\omega(G) = 2$ and $\alpha(G) = 4$.
Interestingly, $R_{3, 5}$ is also the smallest graph $G$ such that
$\theta(G) - \alpha(G) \geq 3$ (see~\cite{GyST}), but we do not use
this fact here.  Observe that $G_{5, 8}$ is obtained from $R_{3, 5}$
by subdividing one edge twice, so that $\omega(G_{5, 8}) = 2$, and
$\theta(G) \geq \lceil |V(G_{5, 8})| / \omega(G_{5, 8})\rceil = 8$.  A
colour class in $G_{5,8}$ is a stable set of size~5 and it is easy to
check that a stable set of size at least~6 in $G_{5,8}$ would contain
a stable set of size 5 of $R_{3, 5}$, a contradiction, so
$\alpha(G_{5, 8}) = 5$.

We now show how to construct a graph $G$
with $\alpha(G) = x$ and $\theta(G) = \lfloor
\frac{8}{5}\alpha(G) \rfloor$ for each integer $x \ge 0$.  Define a graph $G$ consisting of
$k=\lfloor{x\over 5}\rfloor$ disjoint copies of $G_{5,8}$.  So,
$\theta(G)=8k$, $\alpha(G)=5k$.  If $x$ is a multiple of $5$, then
$\theta(G) = \lfloor \frac{8}{5}\alpha(G) \rfloor$ since $\alpha(G)=x$
and $\theta (G)={8\over 5}x$.  If $x\equiv 1$ mod $5$, then add an
isolated vertex to $G$, so that $x=\alpha (G)=5k+1$ and
$\theta(G)=8k+1=\lfloor{8x\over 5}\rfloor$.  If $x\equiv 2$ mod $5$,
then add to $G$ a pentagon, so that $x=\alpha (G)=5k+2$ and $\theta
(G)=8k+3=\lfloor{8x\over 5}\rfloor$.  If $x\equiv 3$ mod $5$, then add
to $G$ a pentagon and an isolated vertex, so that $x=\alpha (G)=5k+3$
and $\theta (G)=8k+4=\lfloor{8x\over 5}\rfloor$.  If $x\equiv 4$ mod
$5$, then add to $G$ two pentagons, so that $x=\alpha(G)=5k+4$ and
$\theta (G)=8k+6=\lfloor{8x\over 5}\rfloor$.

\medskip

To prove the upper bound on $\theta$, we need the next results.  The
following avoids checking small cases. 

\begin{lemma}[Gy\'arf\'as, Seb\H o and Trotignon \cite{GyST}]
  \label{l:gap912}
  If $G$ is a graph on at most 9 vertices, then $\theta(G) - \alpha(G)
  \leq 1$.  If $G$ is a graph on at most 12 vertices, then $\theta(G)
  - \alpha(G) \leq 2$.
\end{lemma}

A graph $G$ is \emph{$\theta$-critical} if for every vertex $v$ of $G$,
$\theta(G-v) < \theta(G)$.   A short proof of the following can be
found in \cite{stehlik:03}.

\begin{theorem}[Gallai \cite{gallai:colorCritical}]
  \label{th:Gcritcial}
  If $G$ is connected and $\theta$-critical, then $$\theta (G) \leq
  \frac{|V(G)| +1}{2}.$$
\end{theorem}

It remains to prove that every graph $G$ in $\cal C$ satisfies
$\theta(G) \leq {8 \over 5} \alpha(G)$.  We prove this by induction on
$|V(G)|$.  For graphs on at most one vertex, the outcome clearly
holds.  If $G$ is not $\theta$-critical, then for some vertex $v$, by
the induction hypothesis, we have $$ \theta(G) = \theta(G-v) \leq {8
  \over 5} \alpha(G-v) \leq {8 \over 5} \alpha(G). $$ If $G$ is
disconnected, then $G$ is the disjoint union of two non-empty graphs
$H_1$ and $H_2$, so by the induction hypothesis $$ \theta (G) =
\theta(H_1) + \theta(H_2) \leq {8 \over 5} (\alpha(H_1) + \alpha(H_2)
) = {8 \over 5} \alpha(G). $$ So we may assume that $G$ is
$\theta$-critical and connected.  By
Theorem~\ref{th:Gcritcial}, $$\theta (G) \leq \frac{|V(G)| +1}{2}.$$
If $\alpha(G) \geq 5$, then $$\theta (G) \leq \frac{|V(G)| +1}{2} \leq
\frac{3\alpha(G) +1}{2} \leq \frac{3\alpha(G) +1}{2} + \frac{\alpha(G)
  -5}{10} = \frac{8}{5} \alpha(G).$$ If $\alpha(G) = 4$, then $|V(G)|
\leq 12$, so by Lemma~\ref{l:gap912}, $\theta(G) \leq \alpha(G) + 2 =
6$, and $\theta(G) \leq 6 < 8\times 4 / 5 = {8 \over 5} \alpha(G) $
is clear.  If $\alpha(G) = 3$, then $|V(G)| \leq 9$, so by
Lemma~\ref{l:gap912}, $\theta(G) \leq \alpha(G) + 1 = 4$, and
$\theta(G) \leq 4 < 8 \times 3 / 5 = {8 \over 5} \alpha(G)$ is
clear.  If $\alpha(G) = 2$, then $|V(G)| \leq 6$, so by
Lemma~\ref{l:gap912}, $\theta(G) \leq \alpha(G) + 1 = 3$, and
$\theta(G) \leq 3 < 8 \times 2 / 5 = {8 \over 5} \alpha(G) $ is
clear.  If $\alpha(G) \leq 1$, then $G$ is a clique, so the outcomes
holds.

\section{A complementary bounding function for eventually identity functions}
\label{sec:evc}

Let $F_c$ denote the class of those $\Nset\rightarrow \Nset$ functions
such that $f(x)=x$ for $x\geq c$.  The following was stated without proof
in \cite{Gyarfas}.

\begin{theorem}\label{old}
  For all $c$ and $g\in F_c$, $g$ is complementary-bounded.
\end{theorem}

\proof We prove by induction on $c$. For $c=1$ only $g(x)=x$
is in $F_c$ and the Perfect Graph Theorem \cite{Lovasz} implies
that $g(x)=x$ is a complementary bounding function.

Suppose that for some $c\ge 1$ every $f\in F_c$ has a complementary
bounding function $f_c^*$ and let $G$ be a graph with $\chi$-bounding
function $g\in F_{c+1}$. Consider a subgraph $H\subseteq G$ with
$\alpha(H)=k$. Let $S=\{s_1,s_2,\dots,s_k\}$ be a stable set in $H$.
Partition $V(H)\setminus S$ into $A_1=N(s_1)$ and for $i=2,\dots ,k$,
$A_i=N(s_i)\setminus (\cup_{j=1}^{i-1} A_j)$.

We claim that for $1\le i \le k$, each $A_i$ induces a graph
$H_i\subset H$ such that $H_i$ has a $\chi$-bounding function in $F_c$. Indeed,
if $H_i$ has a subgraph $G_i$ with $p=\omega(G_i)<\chi(G_i)$ for $p>c$
then $c+1<p+1=\omega(G_i\cup s_i)<\chi(G_i\cup s_i)$ contradicting the
assumption that $G$ is $\chi$-bounded by $g$.  Thus the claim is true so
using $\alpha(H_i)\le k-i+1$ and the induction hypothesis,
$\theta(H_i\cup s_i)=\theta(H_i)\le f_c^*(k-i+1)$, $G$ has a
clique cover with $\sum_{j=1}^k f_c^*(k-j+1)$ cliques.  Thus
$\sum_{j=1}^k f_c^*(k-j+1)$ is a complementary bounding function for~$g$.\qed

For $c=2$, this proof provides the following.

\begin{theorem}\label{upperbf2} $f^*(x) = {x+1\choose 2}$
  is a complementary bounding function for any $f\in F_2$.
\end{theorem}

In fact, the bound provided by Theorem \ref{upperbf2} is at most a
logarithmic factor apart from best possible, since there are
triangle-free graphs $G$ with at least ${c\alpha(G)^2\over
\log{\alpha(G)}}$ vertices, see~\cite{KIM}. For such $G$, $\theta(G)\ge
{c\alpha(G)^2\over 2\log{\alpha(G)}}$.

Theorem \ref{old} invites another question, that of finding $f^*$ for
$f\in F_t$ (with better bounds than Theorem \ref{old}). The first case
beyond the Perfect Graph Theorem (Problem 6.6 in \cite{Gyarfas})
is not even known.

\begin{problem}\label{firststep}
Find $g^*$ for the almost identity function
$$
g(x)=\left\{ \begin{array}{ll}
  3 &\mbox{for } x=2 \\
  x &\mbox{for } x>2\end{array}\right.
$$
\end{problem}

From Corollary \ref{upperbf2}, $g^*(x)\le {x+1\choose 2}$ but it is
possible that $g^*$ is linear. In (\cite{Gyarfas} p. 439), the
conjecture `perhaps $g^*(x)=\lfloor {3x\over 2}\rfloor$ is the truth'
was risked, based on the example of disjoint circuits of length $5$.
This is in fact disproved by the graph
$G_{5, 8}$ represented in Fig.~\ref{fig:g58}.  From
Theorem~\ref{upperbf2} and $G_{5, 8}$ we have $g^*(2)=3$.
While $g^*(3)\ge 4$ is obvious (from the pentagon and an isolated
vertex), it is not clear whether $g^*(3) = 4$.  Theorem~\ref{upperbf2}
suggests the following.

\begin{conjecture}\label{conj:eightfifth} $g^*(x)=\lfloor {8 \over 5} x \rfloor $.
\end{conjecture}

\section{Functions that are not complementary-bounded}
\label{sec:nontb}

In this section, we show that $f(x)=x/\text{polylog}(x)$ is not
complementary-bounded.  We prove this by exhibiting a class of graphs
$\chi$-bounded by $f$ but not $\theta$-bounded.  This family consists
of Schrijver graphs which we define below (in fact, for convenience we
work in the complement, so our graphs will be $\theta$-bounded and not
$\chi$-bounded, but this is clearly equivalent up to a
complementation).

We provide in Lemma~\ref{l:new_Gyarfas} a tool to determine a
$\theta$-bounding function of any graph with ``high'' stability
ratio (ratio between the stability number and number of vertices) and
such that this property is closed under taking induced subgraphs.  It
relies on the following theorem due to Gallai.  We denote by $\nu(G)$
the size of a maximum matching in $G$.

\begin{theorem}[Gallai \cite{gallai:factorCritical}]
  \label{th:gallai}
  If $G$ is a connected graph such that for all vertices $v$, $\nu (G \setminus
  v) = \nu(G)$, then $G$ is factor-critical.
\end{theorem}

\begin{lemma}
  \label{l:new_Gyarfas}
 For every graph $G$
  \[\theta(G) \leq \alpha(G) + \max_{H \subseteq G} (|V(H)| - 2 \alpha(H)).\]
\end{lemma}

\proof We prove the result by induction on $|V(G)|$. It clearly holds
when $|V(G)| \leq 1$.  Note that $\max_{H \subseteq G} (|V(H)| - 2
\alpha(H)) \geq 0$ as we can choose $H$ to be the empty graph.

\noindent{\bf Case 1:} $G$ contains a triangle $T$.

\begin{eqnarray*}
  \theta(G)
  &\leq& 1 + \theta(G\setminus T) \\
  &\leq& 1 + \alpha(G\setminus T) + \max_{H \subseteq G\setminus T}
  (|V(H)| - 2 \alpha(H))\\
  & \leq& \alpha(G) + (1 + \max_{H \subseteq G\setminus T}
  (|V(H)| - 2 \alpha(H)))\\
  & \leq& \alpha(G) + \max_{H \subseteq G}
  (|V(H)| - 2 \alpha(H)).
\end{eqnarray*}

\noindent{\bf Case 2:} There exists a vertex $v\in V(G)$ such that
$\theta(G\setminus v)=\theta(G)$. By the induction hypothesis,
\[
\theta(G) = \theta(G\setminus v) \leq \alpha(G\setminus v) + \max_{H
  \subseteq G\setminus v} (|V(H)| - 2 \alpha(H))\]\[ \leq \alpha(G) +
\max_{H \subseteq G} (|V(H)| - 2 \alpha(H)).
\]

\noindent{\bf Case 3:} $G$ is triangle-free and for every vertex $v\in
V(G)$, $\theta(G\setminus v)<\theta(G)$.  We suppose that $G$ is
connected for otherwise, we obtain the result by the induction
hypothesis on the connected components of $G$.  Observe that for every
vertex $v\in V(G)$, $\theta(G\setminus v)=\theta(G)-1$.  Since $G$ is
triangle-free,
\[
\theta(G)+\nu(G)=|V(G)| \text{ and } \theta(G\setminus v)+\nu(G\setminus v)=|V(G\setminus v)|.
\]

So, for every vertex $v\in V(G)$, $\nu(G\setminus v) =\nu(G)$. By
Theorem~\ref{th:gallai}, $G$ is factor-critical and thus $\theta(G) = (|V(G)| + 1)/2$.

Also, if $G$ is bipartite, then the result holds by
Theorem~\ref{th:k}.  So, from here on, we suppose that $G$ is not
bipartite. An odd cycle $H$ of minimum length in $G$ is chordless,
because a chord would allow us to construct a smaller odd cycle.  It
follows that $\max_{H \subseteq G} (|V(H)| - 2 \alpha(H)) \geq
1$. Now,

\begin{eqnarray*}
\theta(G)
   &=& \frac{|V(G)| +1}{2}
\\ &\leq& \frac{2\alpha(G) + \max_{H \subseteq G} (|V(H)| - 2 \alpha(H)) +1}{2}
\\ &=& \alpha(G) + \max_{H \subseteq G} (|V(H)| - 2 \alpha(H)) +
\frac{-\max_{H \subseteq G} (|V(H)| - 2 \alpha(H)) +1}{2}
\\ &\leq& \alpha(G) + \max_{H \subseteq G} (|V(H)| - 2 \alpha(H)).
\end{eqnarray*}\qed

We now describe subgraphs of Kneser graphs that are not
$\chi$-bounded.  When $n, k$ are integers,
the \emph{Kneser graph} $K\!G_{n,k}$ is the graph whose vertices are
the subsets of $\{1,2,\dots, 2n+k\}$ that have size~$n$, and such that
two vertices are adjacent if they are disjoint sets.  We need several
properties of Kneser graphs.

\begin{lemma}
\label{l:alpha_kneser}
If $H$ is an induced subgraph of $K\!G_{n, k}$, then
\[
\alpha(H) \geq \frac{n}{2n+k} |V(H)|.
\]
\end{lemma}

\proof Let $N = \left | \{ (i, v) : v \in V(H) \text{ and } i \in v\}
\right | = n |V(H)|$.  Suppose that each integer of $\{1, \dots, 2n+k\}$
is contained in less than $({n}/{2n+k}) |V(H)|$ vertices of $H$.  Then
$N = \sum_{1 \leq i \leq 2n+k} |\{v : i\in v\}| < (2n+k) n/(2n+k) |V(H)| = n|V(H)|
= N$, a contradiction.  Therefore, at least one integer of $\{1,
\dots, 2n+k\}$ is contained in at least $\frac{n}{2n+k} |V(H)|$ vertices of
$H$, that form a stable set of $H$.\qed

\begin{lemma}
  \label{l:thetaS}
  If $G$ be an induced subgraph of $KG_{n, k}$, then \[\theta(G) \leq
  \left(1 + \frac{k}{n}\right) \alpha(G).\]
\end{lemma}

\proof By Lemma~\ref{l:new_Gyarfas},

\[\theta(G) \leq \alpha(G) + \max_{H \subseteq G} (|V(H)| - 2
\alpha(H)).\]

Since by Lemma~\ref{l:alpha_kneser},   $|V(H)| \leq (2+\frac{k}{n})
\alpha(H)$, we have:
\[\theta(G) \leq \alpha(G) + \max_{H \subseteq G}
\frac{k}{n}\alpha(H) \leq \alpha(G) + \frac{k}{n} \alpha(G).\]
\qed

An $n$-element subset $S$ of $\{1, \dots, 2n+k\}$ is \emph{sparse}  if it does not contain two neighbors in
the cyclic ordering of $\{1, \dots, 2n+k\}$.  The \emph{Schrijver graph}
$SG_{n,k}$ is the subgraph of $K\!G_{n,k}$ induced by the sparse sets.

\begin{lemma}
  \label{l:count}
$|V(SG_{n, k})| = {n + k \choose n} + {n+k-1 \choose n-1} = {n+k
  \choose k} + {n+k-1 \choose k}$.
\end{lemma}

\proof By \emph{increasing $n$-tuples of $\{1, \dots ,a\}$}, we mean an
$n$-tuple $(i_1, \dots, i_n)$ such that for all $1\leq j \leq n$, we
have $1 \leq i_j \leq a$ and such that $i_1 < \cdots < i_n$.

Let us fisrt count the sparse subsets of $\{1, \dots, 2n+k\}$ that do not
contain $2n+k$.  They are in one to one correspondence with the
increasing $n$-tuples of $\{1, \dots, n+k\}$.  This is clear by
considering the map $(i_1, \dots i_n) \mapsto (i_1, i_2 + 1, \dots,
i_n + n - 1)$.  Therefore there are ${n+k \choose n}$ sparse subsets
of $\{1, \dots, 2n+k\}$ that do not contain $2n+k$.

Let us now count the sparse subsets of $\{1, \dots, 2n+k\}$ that do
contain $2n+k$.  They are in one to one correspondence with the
increasing $(n-1)$-tuples of $\{1, \dots, n+k-1\}$.  This is clear by
considering the map $(i_1, \dots i_{n-1}) \mapsto (i_1+1, i_2 + 2,
\dots, i_{n-1} + n - 1, 2n+k)$.  Therefore, there are ${n+k-1 \choose
  n-1}$ sparse subsets of $\{1, \dots, 2n+k\}$ that do contain
$2n+k$.

The conclusions of the two paragraphs above sum up to the first equality,
and the second follows from ${a+b \choose a} = {a+b \choose b}$.
\qed

The following is the key property of Kneser and Schrijver graphs.

\begin{theorem}[Lov\'asz~\cite{LovaszKneser}, Schrijver \cite{schrijver:knesser}]
  $\chi(K\!G_{n,k})= \chi(SG_{n, k}) = k+2$.
\end{theorem}

We are now ready to  exhibit functions that are not complementary-bounded.

\begin{theorem}
\label{th:unbound1}
Let $h$ be a non-decreasing function such that for all $k$, $h(n^k) \le n/k$ for sufficiently large $n$. Then $x+x/h(x)$ is not complementary-bounded.
\end{theorem}

\proof It is easy to show that for all $k$, $h(2(n+k)^k) \le
\frac{n}{k}$ for sufficiently large $n$ (say for $n \ge N(k)$).  We
can show this by simply choosing the $N(k) = N'(2k+1)$ where $N'$ is
the threshold needed for $h(N^k) \le N/k$ and using the monotonicity
of $h$.  Now, for all $n \ge N(k)$, $h(2(n+k)^k) \le h(n(n+n)^k) =
h(n^{2k+1}) \le n/k$.

By Lemma~\ref{l:count}, for any $k,n$ and any subgraph $H$ of $SG_{2n+k,n}$, we have
\begin{equation*}
\alpha(H) \leq |V(H)| \leq |V(SG_{2n+k, n})| = {n+k
  \choose k} + {n+k-1 \choose k} \leq 2(n+k)^k.
\end{equation*}

We claim the graphs $\cal{S}=\{SG_{2N(k)+k}|k \in \mathbb{N}\}$ are $\theta$-bounded by $f$ (they are not $\chi$-bounded as they are Schrijver graphs).

Since $h$ is non-decreasing, for a subgraph $H$ of $SG_{2n+k} \in \cal{S}$, 
\[
\frac{1}{h(\alpha(H))} \ge \frac{1}{h(2(n+k)^k)} \ge \frac{k}{n}.
\]

So
\[
f(\alpha(H)) = \alpha(H)\left(1+\frac{1}{h(\alpha(H))}\right) \ge \alpha(H)\left(1+\frac{k}{n}\right) \ge \theta(H)
\]
by Lemma \ref{l:alpha_kneser}. Thus, $\cal{S}$ is $\theta$-bounded by $f$ but not $\chi$-bounded, as required.

\qed

\begin{theorem}
\label{th:unbound2}
The function $f(x)=x+x/\log^j(x)$ is not complementary-bounded for any $j\in\Rset$.
\end{theorem}

\proof
We only need to verify that for all $k$, $\log^j(n^k) \le n/k$ for sufficiently large $n$ (and apply Theorem \ref{th:unbound1}).

Given $k$, choose $N$ large enough so $N \ge \log^{2j+1}(N)$ and $\log(N) \ge k$. Then for any $n \ge N$,
\[
k \log^j(n^k) = k^{j+1} \log^j(n) \le \log^{j+1}(n) \log^j(n) = \log^{2j+1}(n) \le n
\]
as required.

\qed

\end{document}